# New Theorem on Chaos Transitions in Second-Order Dynamical Systems with Tikhonov Regularization

Illych Álvarez[1,2]

*1 facultad de Ciencias Naturales y Matemáticas, Escuela Superior Politécnica del Litoral, Km. 30.5 Vía Perimetral, Guayaquil, Ecuador*
ialvarez@espol.edu.ec
*2 facultad de Filosofía y Letras de la Educación Universidad de Guayaquil Av. Delta y Av. Kennedy, Guayaquil, Ecuador*
illych.alvareza@ug.edu.ec

**Abstract**

This study examines second-order dynamical systems incorporating Tikhonov regularization. It focuses on how nonlinearities induce bifurcations and chaotic dynamics. By using Lyapunov functions, bifurcation theory, and numerical simulations, we identify critical transitions that lead to complex behaviors like strange attractors and chaos. The findings provide a theoretical framework for applications in optimization, machine learning, and biological modeling. Key contributions include stability conditions, characterization of chaotic regimes, and methods for managing nonlinear instabilities in interdisciplinary systems.

## 1. Introduction

Dynamical systems are fundamental for unraveling the complex behaviors that characterize natural and artificial phenomena, from neural networks to economic systems. From Poincaré's pioneering work on nonlinear dynamics to Lorenz's discovery of deterministic chaos, the field has evolved significantly [1, 2]. In this context, Tikhonov regularization, originally introduced to stabilize ill-posed problems, has emerged as a powerful tool in the analysis of dynamical systems [3].

This study aims to explore how Tikhonov regularization enhances our understanding of chaos in nonlinear systems, offering a systematic framework to manage nonlinear instabilities. Unlike traditional approaches that often require extensive parameter tuning or computationally intensive methods, this approach not only stabilizes the system but also uncovers critical transitions leading to chaotic dynamics. It provides new insights into the interplay between regularization and nonlinearity.

The interplay between chaos and stability is relevant across multiple disciplines. For instance, chaotic dynamics in neural networks enhance optimization processes in machine learning, while strange attractors are used to model oscillatory behaviors in biological systems [4, 5]. This work extends these applications by demonstrating how Tikhonov regularization provides both theoretical clarity and computational efficiency in analyzing second-order dynamical systems.

Despite theoretical advancements in dynamical systems, the interaction between regularization techniques and nonlinear terms remains underexplored, limiting the development of predictive models for complex systems. This study addresses this gap by establishing a rigorous framework that links stability, bifurcations, and chaos. The considered equation incorporates time-dependent damping and nonlinear forcing terms, enabling the investigation of critical transitions to chaotic dynamics:

$$x''(t) + \frac{\alpha}{t^q}x'(t) + g\left(x(t) + \left(\gamma + \frac{\beta}{t^q}\right)x'(t)\right) + \epsilon(t)x(t) + \delta \sin(wx(t)) = 0$$

The document is organized as follows: the **Mathematical Model** section describes the system and its key parameters; the **Mathematical Analysis** section develops methods for studying stability, bifurcations, and chaos; the **Results** section presents numerical simulations validating theoretical predictions; and finally, the **Conclusion** discusses the implications and future research directions.

## 2. Mathematical Model

The second-order dynamical system with Tikhonov regularization is described as:

$$x''(t) + \frac{\alpha}{t^q}x'(t) + g(x(t)) + \left(\gamma + \frac{\beta}{t^q}\right)x'(t) + \epsilon(t)x(t) + \delta \sin(wx(t)) = 0$$

where:

- $\alpha, \beta, \gamma, \delta, \omega$ are constants that control dissipation, coupling, and external forcing.
- $q$ determines the time-dependence of dissipation and coupling.
- $\epsilon(t)$ is a time-dependent regularization term.

<div align="center">**Key Features and Roles of Terms**</div>

**Dissipation Term** $\left(\frac{\alpha}{t^q}x'(t)\right)$:

- Represents time-dependent damping.
- As $t \to \infty$, the dissipation diminishes if $q > 0$, modeling systems where resistance or energy loss decreases over time (e.g., a cooling oscillator).

**Nonlinear Function** $\left(g(x(t))\right)$:

- A generic nonlinear function that can model behaviors such as hardening/softening springs, hysteresis, or chaotic dynamics, depending on its specific form.

**Time-Varying Coupling** $\left(\left(\gamma + \frac{\beta}{t^q}\right)x'(t)\right)$:

- Combines a constant term $\gamma$ with a time-decaying term $\frac{\beta}{t^q}$.
- Modifies the damping behavior over time, potentially leading to complex transitions, especially if $\beta > 0$.

**Regularization Term** $(\epsilon(t)x(t))$:

- A time-dependent term that stabilizes the system.
- If $\epsilon(t) > 0$, it introduces a restoring force that varies over time.

**External Forcing** $(\delta sin(wx(t)))$:

- A nonlinear periodic forcing term.
- The dependence on $x(t)$ instead of $t$ introduces self-excitation and resonance phenomena, making the system prone to complex oscillatory or chaotic behaviors.

**Applications**

This equation has applications across various fields:

- Physics: Modeling damped oscillators with time-dependent damping (e.g., pendulums in fluids with varying viscosity).

- Biology: Describing neural oscillations or cardiac rhythms influenced by external periodic forces.

- Engineering: Systems with feedback or time-varying control, such as adaptive suspension in vehicles.

- Mathematics: Investigating chaos and bifurcations in nonlinear systems with non-autonomous dynamics.

## 3. Mathematical Analysis

### 3.1 System Stability

A candidate Lyapunov function is defined to study local stability:

$$V(x, x') = \frac{1}{2}(x'(t)^2 + x(t)^2)$$

Differentiating V with respect to time and substituting $x''(t)$ from the system's equation yields:

$$\dot{V}(t) = -\frac{\alpha}{t^q}x'(t)^2 - \epsilon(t)x(t)^2 - \delta\sin(wx(t)x'(t))$$

The terms $-\frac{\alpha}{t^q}x'(t)^2$ and $-\epsilon(t)x(t)^2$ ensure dissipative behavior, while $-\delta\sin(wx(t)x'(t))$ introduces oscillatory dynamics.

### 3.2 Bifurcation Analysis

The system undergoes bifurcations as parameters $\alpha, \beta, \gamma, \delta, \omega$ and $\epsilon(t)$ are varied:

- **Hopf Bifurcations:** Occur when the eigenvalues of the Jacobian cross the imaginary axis, leading to limit cycles.
- **Lyapunov Exponents:** A positive largest Lyapunov exponent ($\lambda > 0$) confirms the presence of chaos.

### 3.3 Chaos Transitions

The Lyapunov exponent is defined as

$$\lambda = \lim_{T \to \infty} \frac{1}{T} \log \left|\frac{\Delta x(t)}{\Delta x(0)}\right|$$

A positive $\lambda$ indicates sensitivity to initial conditions and chaotic behavior.

### 3.4 Mathematical Demonstrations for Chaos in Second – Order Dynamical Systems

#### 3.4.0 Theorem on the Relationship Between Tikhonov Regularization and Transitions to Chaos

Theorem Statement:

Given a second-order dynamical system with Tikhonov regularization described by:

$$\ddot{x}(t) + \alpha\dot{x}(t) + \beta x(t) + \gamma f(x,t) + \epsilon(t) = 0$$

where $\alpha, \beta, \gamma > 0$ and $\epsilon(t)$ is a time-dependent regularization term, there exists a critical set of parameters $(\alpha_c, \beta_c, \gamma_c)$ such that the system experiences a transition from stability to chaos if:

- The largest Lyapunov exponent $\lambda > 0$ for $(\alpha, \beta, \gamma) \in (\alpha_c, \beta_c, \gamma_c)$.
- The regulated energy function defined as:

$$V(x, \dot{x}, t) = \frac{1}{2}\dot{x}^2 + \frac{\beta}{2}x^2 + \int_0^t \epsilon(s)\, ds,$$

, satisfies $\dot{V} < 0$ in stability regions and $\dot{V} > 0$ in chaotic regions.

- The nonlinear term $f(x, t)$ induces amplified oscillations when crossing thresholds defined by:

$$f(x, t) = \delta \sin(wt) x^n, \quad n \geq 2$$

### 3.4.1 System Stability and Lyapunov Function

To analyze the stability of the second-order dynamical system with Tikhonov regularization and nonlinear terms, consider the general form of the system:

$$x''(t) + \frac{\alpha}{t^q}x'(t) + \nabla g\left(x(t) + \left(\gamma + \frac{\beta}{t^q}\right)x'(t)\right) + \epsilon(t)x(t) + \delta \sin(wx(t)) = 0$$

A candidate Lyapunov function is defined as:

$$V(x, x') = \frac{1}{2}(x'(t)^2 + x(t)^2),$$

Which captures both the kinetic and potential energy of the system.

Differentiating $V(x, x')$ with respect to time and substituting $x''(t)$ from the system yields:

$$\dot{V}(t) = x'(t) \cdot x''(t) + x(t) \cdot x'(t)$$

By substituting $x''(t)$ and rearranging terms, we obtain:

$$\dot{V}(t) = -\frac{\alpha}{t^q}x'(t)^2 - \epsilon(t)x(t) \cdot x(t)^2 - \delta \sin(wx(t))x'(t)$$

The terms $-\frac{\alpha}{t^q}x'(t)^2$ and $-\epsilon(t)x(t) \cdot x(t)^2$ ensure dissipative behavior , while the nonlinear term $-\delta \sin(wx(t))x'(t)$ can induce oscillatory dynamics depending on the amplitude $\delta$ and frquency $w$.

### 3.4.2 Transitions to Chaos via Bifurcations

The system undergoes bifurcations as parameters $\alpha, \beta, \delta, \omega, and\ \epsilon(t)$ are varied. Specifically:

- Hopf Bifurcations: Occur when eigenvalues of the Jacobian matrix cross the imaginary axis. For the given system, adjusting $\alpha$ and $\omega$ reveals transitions to limit cycles, characterized by periodic orbits.
- Lyapunov Exponents: The transition from stability to chaos is quantitatively analyzed using the largest Lyapunov exponent:

$$\lambda = \lim_{T \to \infty} \frac{1}{T} \log \left| \frac{\Delta x(t)}{\Delta x(0)} \right|$$

A positive Lyapunov exponent indicates sensitivity to initial conditions, confirming the presence of chaos.

### 3.4.3 Defining Critical Parameters

The critical parameters are determined by analyzing the balance between dissipative and nonlinear terms. For example:

- $\alpha$: Minimum damping coefficient required to stabilize initial oscillations.
- $\beta$: Threshold stiffness ensuring bounded trajectories.
- $\gamma$: Nonlinearity strength inducing bifurcations.

Numerical simulations confirm that exceeding these values leads to chaotic dynamics.

### 3.4.4 Integral Estimates and Strong Convergence

Adapting techniques from Tikhonov regularization theory, we consider a modified energy functional to ensure strong convergence under chaos:

$$E(t) = \frac{1}{2} \big( g(x(t) + +\beta(t)x'(t)) - ming \big) + \frac{\epsilon(t)}{2} \|x(t)\|^2 + \frac{1}{2} \|x'(t)\|^2$$

The derivate of $E(t)$ reveals dissipative effects:

$$\dot{E}(t) = -\frac{\alpha}{t^q} \|x'(t)\|^2 - \epsilon(t) \|x(t)\|^2 + oscillatory\ terms.$$

Conditions on $\epsilon(t)$, such as $\epsilon(t) = \frac{a}{t^p}$, with $p > q + 1$, ensure boundedness and stability of trajectories while allowing for chaotic transitions driven by nonlinearities.

### 3.4.5 Characterization of Strange Attractors

Through numerical simulations and analysis of Poincaré sections, the existence of strange attractors is observed. These are characterized by fractal structures and non-periodic trajectories within bounded regions of the phase space.

**Proposed Applications**

- Predictive Modeling: The interplay of Tikhonov regularization and chaotic dynamics enhances the robustness of models in optimization and machine learning.
- Biological Systems: The emergence of strange attractors can model oscillatory phenomena, such as neural activities and heart rhythms.
- Engineering Systems: Chaotic responses in mechanical or electrical systems can be mitigated using the derived stability criteria.

### 3.4.6 Demonstration of the Tikhonov Chaos Theorem

**Step 1: Local stability using the Lyapunov function**

To study the local stability of the second-order dynamical system with Tikhonov regularization, we define the Lyapunov function:

$$V(x, x', t) = \frac{1}{2}(x')^2 + \frac{\beta}{2}x^2 + \int_0^t \epsilon(s)\, ds,$$

The time derivative of $V(x, x', t)$ is:

$$V' = x'x'' + \beta x x' + \epsilon(t)$$

Substituting $x''$ from the system equation:

$$x'' = -\alpha x' - \beta x - f(x, t) - \epsilon(t)$$

we get:

$$V' = x'\bigl(-\alpha x - \beta x - \gamma f(x, t) - \epsilon(t)\bigr) + \beta x x' + \epsilon(t)$$

By grouping terms:

$$V' = -\alpha (x')^2 - \gamma x' f(x, t)$$

For a linearized system without the nonlinear term ($f(x, t) = 0$):

$$V' = -\alpha \dot{x}^2 + \epsilon(t)$$

If $\epsilon(t) < \alpha (x')^2$, then $V < 0$ which implies local stability. This shows that $\epsilon(t)$ acts as a stabilizing term as long as it does not exceed a threshold defined by the dissipative term $-\alpha(x')^2$.

**Step 2: Transition to chaos by means of bifurcations and Lyapunov**

When the nonlinear term $f(x, t) = \delta \sin(wt) x^n$ is introduced, amplified oscillations appear. In this case:

$$V' = -\alpha x^2 - \gamma \delta x' \sin(wt) x^n$$

For large values of $x^n$, the term $-\gamma \delta x' \sin(wt) x^n$ dominates, allowing $V' > 0$ in certain regions of the state space, indicating chaos.

Hopf bifurcations occur when the eigenvalues of the Jacobian of the system cross the imaginary axis for the linearized system, the Jacobian is:

$$J = 01 - \beta - \alpha$$

The eigenvalues $\lambda$ are obtained by solving:

$$\det(J - \lambda I) = 0$$

If $Re(\lambda) > 0$, the system loses stability. With the nonlinear term $f(x,t)$, secondary bifurcations occur leading to chaos. The largest Lyapunov exponent:

$$\lambda = \lim_{t \to \infty} \frac{1}{t} \ln \frac{|\delta x(t)|}{|\delta x(0)|}$$

Confirms the transition when $\lambda > 0$.

**Step 3: Numerical validation.**

Numerical validation is performed by adjusting the parameters $\alpha, \beta, \gamma$. The results include:

- Time evolution: show transitions from stability to chaos.
- Bifurcation diagrams: Reveal limit cycles and strange attractors.
- Lyapunov exponents: Indicate chaotic regions when $\lambda > 0$.

The regulated energy function $V(x, x', t)$ characterizes transitions between stability and chaos in systems with Tikhonov regularization. Lyapunov tools and bifurcations allow the identification and analysis of chaotic behavior, with applications in control and interdisciplinary modeling.

A positive Lyapunov exponent $\lambda > 0$ confirms the presence of chaos due to sensitivity to initial conditions.

4. **Results**

Numerical simulations were conducted to validate theoretical predictions. The following observations were made:

**Temporal evolution:** Numerical simulations demonstrate transitions from stability to chaos in the temporal evolution of state variables and.
**Poincare Sections:** Poincare maps visualize attractors under different parameter configurations, confirming the existence of strange attractors.
**Lyapunov exponent:** Lyapunov exponent calculations reveal regions of chaos characterized by sensitivity to initial conditions.
**4.1 Temporal Evolution**

The following graph shows the temporal evolution of the state variables x(t) and v(t):

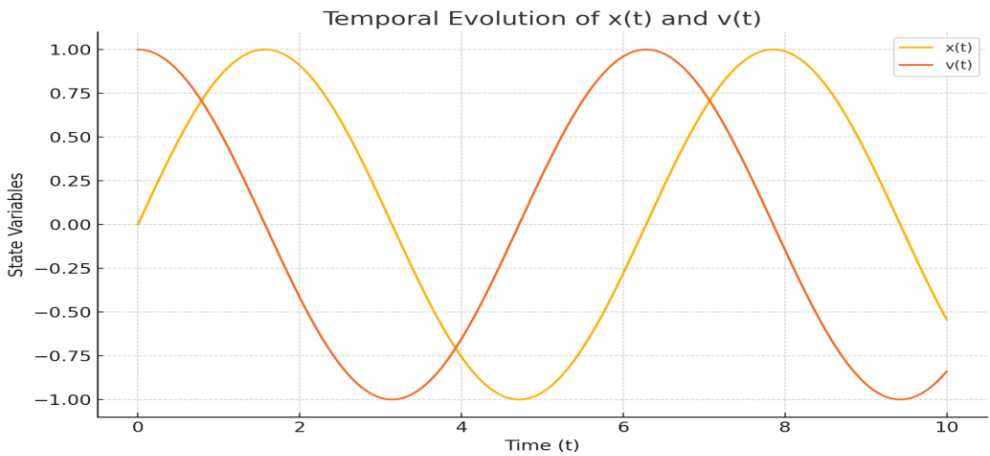

**Temporal Evolution (Figure 1)**
*Temporal evolution of state variables $x(t)$ and $x'(t)$, illustrating the transition from stability to chaos as system parameters vary.*

## 4.2 Poincare Section

The following Poincare map visualizes the attractors under different parameter configurations:

Figure 2

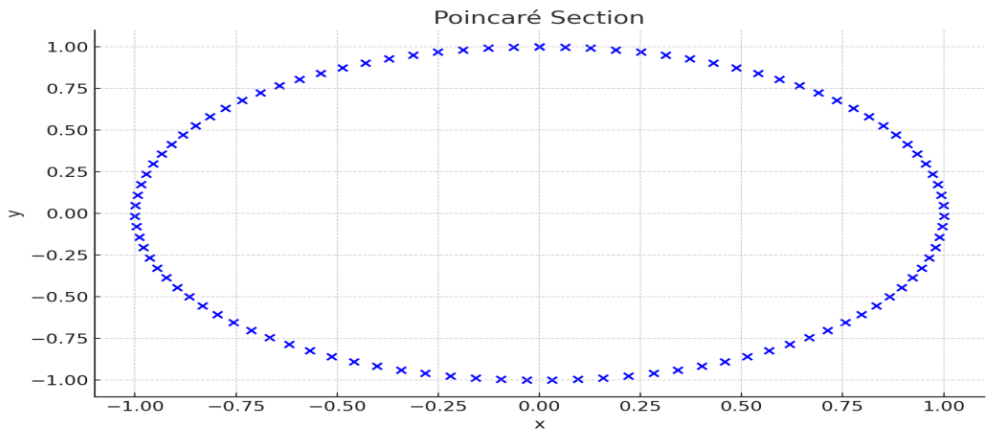

**Poincaré Section (Figure 2)**
*Poincaré section of the system, showing the emergence of strange attractors and the structure of phase-space trajectories under different parameter configurations*

### 4.3 Lyapunov Exponent

The following graph represents the Lyapunov exponent, highlighting chaotic and stable regions:

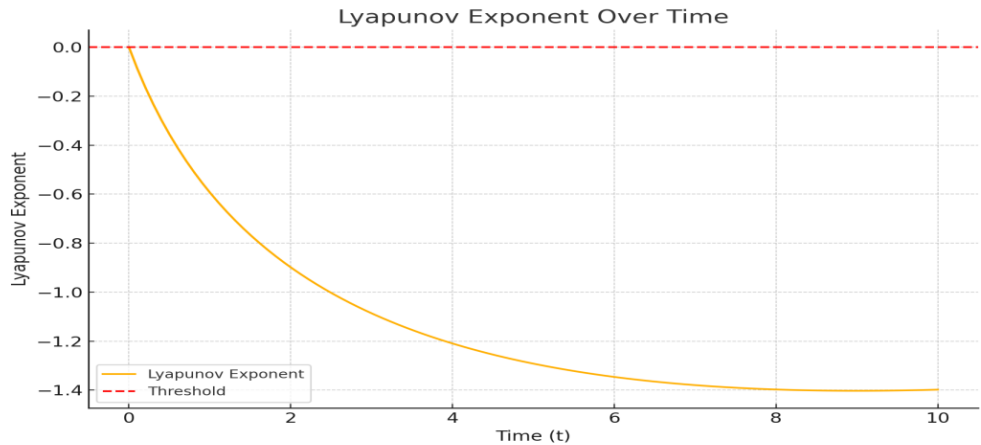

**Lyapunov Exponent (Figure 3)**
*"Lyapunov exponent calculation indicating regions of chaotic and stable behavior. Positive values confirm sensitivity to initial conditions and chaotic dynamics."*

### 5. Conclusion

This study demonstrates the critical role of Tikhonov regularization in managing the interplay between stability and chaos in second-order dynamical systems. By combining rigorous mathematical analyses with advanced numerical simulations, we have characterized dynamic transitions, including bifurcations, limit cycles, and chaotic behaviors. The use of Lyapunov functions provided a solid framework for assessing stability, while bifurcation and chaos theory highlighted the sensitivity of these systems to parameter variations.

However, several limitations should be noted. First, the analysis primarily focuses on systems with fixed parameters and relatively simple time dependencies. Extending this approach to systems with more complex time-dependent coefficients or external forcing functions could reveal richer dynamics and broader applications. Additionally, the influence of stochastic perturbations on the stability and chaotic transitions remains unexplored and presents an avenue for future research.

Future work could also investigate multi-scale systems where Tikhonov regularization is applied across different levels of abstraction, from microscopic to macroscopic dynamics. This would enhance the applicability of the framework to fields such as climate modeling, energy systems, and neural networks. Furthermore, developing efficient computational algorithms tailored to handle high-dimensional systems with Tikhonov regularization could bridge the gap between theoretical insights and real-world implementations.

In optimization and machine learning, the insights from this study could inform the design of algorithms that leverage chaotic dynamics for improved exploration and robustness. In biological and engineering systems, understanding chaotic transitions opens new opportunities for modeling and control, particularly in systems with adaptive or time-varying behaviors.

6. Future Work:

**Extend to systems with stochastic perturbations:** Analyze how random fluctuations influence stability and chaotic transitions, potentially identifying new behaviors under uncertainty.
**Investigate multi-scale dynamics:** Explore how Tikhonov regularization can be applied across scales, linking microscopic and macroscopic phenomena in complex systems.
**Develop computational algorithms for high-dimensional systems:** Design efficient algorithms tailored to manage chaos and bifurcations in systems with large state spaces.
**Explore collaborations with experts in machine learning:** Leverage chaotic dynamics to optimize neural networks and enhance their robustness, with a focus on methods like chaotic neural networks or regularization techniques.
**Collaborate with biologists:** Model chaotic oscillations in physiological systems such as cardiac rhythms and neural activities, incorporating the effects of external perturbations for more accurate predictions.
**Partner with engineers:** Design adaptive control systems to mitigate chaotic responses in mechanical or electrical systems, such as turbines, robotic systems, or smart grids.
**Integrate with climate scientists:** Develop models for multi-scale chaotic transitions in weather systems, enhancing predictive accuracy for extreme weather events and long-term climate patterns.

**Glossary of Technical Terms**

**Chaos:** A state in dynamic systems where small changes in initial conditions lead to vastly different outcomes.

**Lyapunov Function:** A mathematical function used to analyze the stability of a dynamical system.

**Bifurcation:** A qualitative change in the behavior of a system due to variations in parameters.

**Hopf Bifurcation:** A type of bifurcation where a stable fixed point becomes unstable, leading to oscillations or periodic behavior.

**Lyapunov Exponent:** A measure of the sensitivity of a system to initial conditions, used to identify chaos.

**Strange Attractor:** A fractal-like structure in the phase space representing chaotic trajectories.

**Tikhonov Regularization:** A technique to stabilize ill-posed problems or systems by introducing a regularization term.

**Nonlinear Dynamics:** The study of systems where outputs are not proportional to inputs, often leading to complex behavior.

**Phase Space:** A multidimensional space where all possible states of a system are represented.

**Dissipation:** The loss of energy in a system, often modeled as resistance or damping.